\documentclass{article}
\title{A characterization of the category ${\bf {\it Q}}$-${\bf TOP}$}
\author{Sheo Kumar Singh\thanks{sheomathbhu@gmail.com}\\ {\it Department of Mathematics},\\{\it Banaras Hindu University},\\{\it Varanasi-221005, India}\\\\
Arun K. Srivastava\thanks{arunksrivastava@gmail.com}\\
{\it Department of Mathematics}\\{\it and}\\{\it Centre for Interdisciplinary Mathematical Sciences},\\{\it Banaras Hindu University},\\{\it Varanasi-221005, India} 
 }
\date{}
\usepackage{amsmath,amssymb,mathrsfs}
\newtheorem{thm}{Theorem}[section]

\newtheorem{rem}{Remark}[section]
\newtheorem{dfn}{Definition}[section]

\usepackage[all]{xy} 
\begin{document}
\maketitle
\section{Introduction} E.G. Manes in \cite{mane} gave a somewhat `axiomatic' characterization (upto isomorphism) of the category ${\bf TOP}$ of topological spaces, among a certain class of categories, satisfying some conditions, in which a Sierpinski space-like object, played a key role. Following that, the category $[0, 1]$-${\bf TOP}$ of $[0, 1]$-topological spaces (known more commonly as `fuzzy' topological spaces) was analogously characterized by Srivastava \cite{sriv}. S.A. Solovyov \cite{solo} has recently introduced the notion of $Q$-topological spaces ($Q$ being a fixed member of a variety of $\Omega$-algebras). In this note, we give a characterization of the category ${\bf {\it Q}}$-${\bf TOP}$ of $Q$-topological spaces in which the so-called $Q$-Sierpinski space, introduced in \cite{solo}, plays a key role.

\section{The category ${\bf {\it Q}}$-${\bf TOP}$}
We begin by recalling the (well-known) notions of $\Omega$-algebras and their homomorphisms. For details, see \cite{jaco}.\\
\newpage
\begin{dfn} {\rm(\cite{solo})}
\begin{itemize}
\item Let $\Omega = (n_\lambda)_{\lambda\in I}$ be a class of cardinal numbers. An {\bf $\Omega$-algebra} is a pair $(A, (\omega _{\lambda}^{A})_{\lambda\in I})$ consisting of a set $A$ and a family of maps $\omega _{\lambda}^{A}: A^{n_\lambda}\rightarrow A$. Any $B(\neq \phi)\subseteq A$ is called a {\bf subalgebra} of $(A, (\omega _{\lambda}^{A})_{\lambda\in I})$ if for any $\lambda \in I$ and any $(b_i)_{i \in n_\lambda}\in B^{n_\lambda}$, $\omega _{\lambda}^{A}((b_i)_{i \in n_\lambda})\in B$. Given $S\subseteq A$, the intersection of all the subalgebras of $(A, (\omega _{\lambda}^{A})_{\lambda\in I})$ containing $S$ is clearly a subalgebra of $(A, (\omega _{\lambda}^{A})_{\lambda\in I})$. We shall denote it as $<S>$.  

\item Given $(A, (\omega _{\lambda}^{A})_{\lambda\in I})$ and $(B, (\omega _{\lambda}^{B})_{\lambda\in I})$, a map $f: A\rightarrow B$ is called an {\bf $\Omega$-algebra homomorphism} if for every $\lambda \in I$ the following diagram commutes: 
$$\xymatrix{A^{n_\lambda}\ar[d]_{\omega _{\lambda}^{A}}\ar[r]^{f^{n_\lambda}}&B^{n_\lambda}\ar[d]^{\omega _{\lambda}^{B}}\\
            A\ar[r]_{f}&B}$$.\\
 Let ${\bf Alg}(\Omega)$ denote the category of $\Omega$-algebras and $\Omega$-algebra homomorphisms.
            
\item A {\bf variety} of $\Omega$-algebras is a full subcategory of ${\bf Alg}(\Omega)$ which is closed under the formation of products \footnote{The category ${\bf Alg}(\Omega)$ is closed under products}, subalgebras and homomorphic images.\\
  
  \quad  {\rm {\bf Throughout, $Q$ denotes a fixed member of a fixed variety of $\Omega$-algebras}}.

\item Given a set $X$, a subset $\tau$ of $Q^X$ is called a {\bf $Q$-topology} on $X$ if $\tau$ is a subalgebra of $Q^X$; in such a case, the pair $(X, \tau)$ is called a {\bf $Q$-topological space}. 

\item Given two $Q$-topological spaces $(X, \tau)$ and $(Y, \eta)$, a {\bf $Q$-continuous function} from $(X, \tau)$ to $(Y, \eta)$ is a function $f: X\rightarrow Y$ such that $f^{\leftarrow}(\alpha) \in \tau, \forall \alpha \in \eta$, where $f^{\leftarrow}(\alpha)= \alpha \circ f$.

\end{itemize}
\end{dfn}

\quad It is evident that all $Q$-topological spaces, together with $Q$-continuous maps, form a category, which we shall denote as ${\bf {\it Q}}$-${\bf TOP}$.\\

\quad The following categorical concepts are from Manes \cite{mane}
\begin{dfn} {\rm(\cite{mane})}
\begin{enumerate}
\item A {\bf category} $C$ {\bf of sets with structure} is defined through the following descriptions of its objects {\rm (}`$C$-structured sets'{\rm )} and morphisms {\rm (}`$C$-admissible maps'{\rm )}, satisfying the two axioms given below:\\
\indent A class $C(X)$ of `$C$-structures' is assigned with each set $X$ and a `$C$-structured set' is a pair $(X, s)$, with $s\in C(X)$.\\
\indent A subset $C(s, t)$ of the set of all functions from $X$ to $Y$ is assigned with each pair of $C$-structured sets $(X, s)$ and $(Y, t)$ and a `$C$-admissible map' from $(X, s)$ to $(Y, t)$ is any $f\in C(s, t)$ {\rm (}in which case we write ``$f:(X, s)\rightarrow (Y, t)$"{\rm )}.\\
\indent The axioms are:\\
{\rm Axiom $A_1$}: If $f:(X, s)\rightarrow (Y, t)$ and $g:(Y, t)\rightarrow (Z, u)$ then also $g\circ f:(X, s)\rightarrow (Z, u)$.\\
{\rm Axiom $A_2$}: Given a bijection $f:X\rightarrow Y$ and $t\in C(Y)$, there exists a unique $s\in C(X)$ such that $f:(X, s)\rightarrow (Y, t)$ and $f^{-1}:(Y, t)\rightarrow (X, s)$.
\item Given a category ${\bf C}$ of structures {\rm (}as defined above{\rm )}, a family ${\cal F}=\{f_j:(X, s)\rightarrow (Y_j, t_j)| j\in J\}$ of $C$-admissible maps is said to be {\bf optimal} if for each $C$-structured set $(Z, u)$ and a function $g:Z\rightarrow X$,  $g:(Z, u)\rightarrow (X, s)$ iff $f_j\circ g:(Z, u)\rightarrow (X, s), \forall j\in J$. Further, if for a family ${\cal F}=\{f_j:X\rightarrow (Y_j, t_j)| j\in J\}$ of functions, where $X$ is a set and each $(Y_j, t_j)$ is a $C$-structured set, there exists $s\in C(X)$ such that the family $\{f_j:(X, s)\rightarrow (Y_j, t_j)| j\in J\}$ is optimal, then $s$ is called an {\bf optimal lift} of the family ${\cal F}$.

\item An object $S= (S, u)$ in a category $C$ of sets with structure is called a {\bf Sierpinski object} if for every $C$-object $X=(X, s)$, the family of all $C$-admissible maps from $X$ to $S$ is optimal.

\end{enumerate}
\end{dfn}

\begin{rem}
It is easy to verify that, in ${\bf {\it Q}}$-${\bf TOP}$, the optimal lift of a family ${\cal F}=\{f_j:X\rightarrow (Y_j, t_j)| j\in J\}$ of functions, where $X$ is a set and each $(Y_j, t_j)$ is a ${\bf {\it Q}}$-${\bf TOP}$-object, is precisely the smallest $Q$-topology on $X$ making each $f_i$ $Q$-continuous.
\end{rem}

\qquad Both ${\bf TOP}$ and ${\bf [0, 1]}$-${\bf TOP}$ are categories of sets with structures and the usual two-point Sierpinski space and the fuzzy Sierpinski space (of \cite{sriv}) are Sierpinski objects in these categories.

\qquad We verify (on expected lines) that ${\bf {\it Q}}$-${\bf TOP}$ is also a category of sets with structures. For each set $X$, $C(X)$ is the family of all $Q$-topologies on $X$ and the $C$-admissible maps are just the $Q$-continuous maps. Given a bijection $f:X\rightarrow Y$ and $t\in C(Y)$, there exists a unique $s\in C(X)$, namely $s= <\mathscr{A}>$, with $\mathscr{A}= \{q\circ f | q\in t\}$, such that both $f:(X, s)\rightarrow (Y, t)$ and $f^{-1}:(Y, t)\rightarrow (X, s)$ are $Q$-continuous. Also, every family ${\cal F}=\{f_j:X \rightarrow (Y_j, t_j)| j\in J\}$ of functions, where $X$ is a set and $(Y_j, t_j)$, $j\in J$ are the $Q$-topological spaces, has an optimal lift $s$, namely $s= <\mathscr{S}>$, with $\mathscr{S}= \{q\circ f_{j} | q\in t_{j}, j\in J\}$. The {\it product} of a family $\{(X_i,t_i) | i\in I\}$ of $Q$-topological spaces is the $Q$-topological space $(\prod X_i, t)$, where $t$ is the optimal lift of the family of all the projection maps $p_i: \prod X_i \rightarrow (X_i, t_i)$. In fact, $(\prod X_i, t)$ is the categorical product of the ${\bf {\it Q}}$-${\bf TOP}$-objects $(X_i, t_i)$, $i\in I$.

\section{$Q$-Sierpinski Space}
Let $S=Q$. It is obvious that $u= <id>$, where $id \in Q^Q$ is the identity function, is a $Q$-topology on $S$. $(S, u)$ has been called the ${\bf {\it Q}}$-{\bf Sierpinski space} in \cite{solo}. Theorem 3.2 establishes the appropriateness of this concept of $Q$-Sierpinski space. But first, we state the following easy-to-verify result.

\begin{thm}
For any $Q$-topological space $(X, \tau)$, $p\in \tau$ iff $p: (X, \tau)\rightarrow (S, u)$ is $Q$-continuous.
\end{thm}

\begin{thm}
$(S, u)$ is a Sierpinski object in ${\bf {\it Q}}$-${\bf TOP}$.
\end{thm}

{\bf Proof:} Let $(X, \tau)$ be a ${\bf {\it Q}}$-${\bf TOP}$-object and ${\cal F}=\{f:(X, \tau)\rightarrow (S, u) | f$ is $Q$-continuous$\}$ be the family of all $Q$-continuous maps from $(X, \tau)$ to $(S, u)$. To show that this family is optimal, let $\tau'$ be any other $Q$-topology on $X$ making each $f\in {\cal F}$ $Q$-continuous. If $p\in \tau$ then, as $id \in u$, $p:(X, \tau)\rightarrow (S, u)$ must be $Q$-continuous and so $p:(X, \tau')\rightarrow (S, u)$ is also $Q$-continuous. But then $p\in \tau'$. Thus $\tau \subseteq \tau'$, showing that $\tau$ is the smallest $Q$-topology on $X$ making each $f\in {\cal F}$ $Q$-continuous. $\square$\\


\section{A characterization of ${\bf {\it Q}}$-${\bf TOP}$}

\begin{thm}
A category ${\bf C}$ of sets with structures is isomorphic to the category ${\bf {\it Q}}$-${\bf TOP}$ if and only if ${\bf C}$ contains an object $(S, u)$, with $S=Q$, satisfying the following conditions:
\begin{enumerate}
\item every family $f_i:X\longrightarrow (S, u)$ has an optimal lift,
\item the maps $\omega_\lambda:(S, u)^{n_\lambda} \longrightarrow (S, u)$ are ${\bf C}$-admissible for each $\lambda \in I$, where $(S, u)^{n_\lambda}= (S^{n_\lambda}, u_{n_\lambda})$, with $u_{n_\lambda}$ being the optimal lift of all the projection maps from $S^{n_\lambda}$ to $(S, u)$,
\item $(S, u)$ is a Sierpinski object in ${\bf C}.$
\end{enumerate}
\end{thm}

{\bf Proof:} First, it is clear that ${\bf {\it Q}}$-${\bf TOP}$ satisfies (1)
 and (3). That ${\bf {\it Q}}$-${\bf TOP}$ also satisfies $(2)$ is shown as follows. Note that $u_{n_\lambda}$, the product $Q$-topology on $S^{n_\lambda}$ has a generating set  consisting precisely of all the projection maps from $S^{n_\lambda}$ to $S$. To show that $\omega_\lambda:(S, u)^{n_\lambda} \longrightarrow (S, u)$ is $Q$-continuous, it suffices to show that $\omega_\lambda \in u_{n_\lambda}$. Now if $a=(a_i)_{i\in n_\lambda}\in S^{n_\lambda}$, then ${\omega_\lambda}^{S^{S^{n_\lambda}}} ((p_i)_{i\in n_\lambda})(a)= \omega_\lambda ^S((p_i(a))_{i\in n_\lambda})= \omega_\lambda ^S((a_i)_{i\in n_\lambda})= \omega_\lambda (a)$, where ${p_i}'s$ are the projection maps from $S^{n_\lambda}$ to $S$. Hence $\omega_\lambda ={\omega_\lambda}^{S^{S^{n_\lambda}}}((p_i)_{i\in n_\lambda})\in u_{n_\lambda}$. Thus $\omega_\lambda$ is $Q$-continuous, $\forall \lambda \in I$.

For each set $X$, let $Q(X)$ denote the set of all $Q$-topologies on $X$. It will suffice to prove that, for each set $X$, there exists a bijection $\Phi_X: C(X)\rightarrow Q(X)$ such that $f: (X, s)\rightarrow (Y, t)$ is $C$-admissible iff $f: (X, {\Phi_X}(s))\rightarrow (Y, {\Phi_Y}(t))$ is $Q$-continuous, i.e., iff $q\circ f\in {\Phi_X}(s)$ for each $q\in {\Phi_Y}(t)$. For each set $X$ and $s\in C(X)$, put ${\Phi_X}(s)= \{p\in Q^X | p:(X, s)\rightarrow (S, u)$ is $C$-admissible$\}$. We show that ${\Phi_X}(s)$ is a $Q$-topology on $X$, i.e., ${\Phi_X}(s)$ is a subalgebra of $Q^X$, or that ${\Phi_X}(s)\in Q(X)$. For a given $\lambda \in I$, consider $(p_i)_{i\in {n_\lambda}}\in (Q^X)^{n_\lambda}$, with $p_i \in {\Phi_X}(s)$, $\forall i\in {n_\lambda}$. We wish to show that $\omega_\lambda ^{Q^X}((p_i)_{i\in n_\lambda})\in {\Phi_X}(s)$, i.e., $\omega_\lambda ^{Q^X}((p_i)_{i\in n_\lambda})$ is also a $C$-admissible map from $(X, s)$ to $(S, u)$. The map $f:(X, s)\rightarrow (S, u)^{n_\lambda}$, defined by $f(x)(j)= p_j(x)$, for each $x\in X$ and $1 \leq j\leq {n_\lambda}$, is easily seen to be $C$-admissible. Here we note that $\omega_\lambda ^{Q^X}((p_i)_{i\in n_\lambda})=\omega_\lambda \circ f$, as $(\omega_\lambda \circ f)(x)=\omega_\lambda(f(x))= \omega_\lambda ^Q((f(x)(i))_{i\in n_\lambda})= \omega_\lambda ^{Q}((p_i(x))_{i\in n_\lambda})=\omega_\lambda ^{Q^X}((p_i)_{i\in n_\lambda})(x)$. So, $\omega_\lambda ^{Q^X}((p_i)_{i\in n_\lambda})$ is also a $C$-admissible, being the composition of two $C$-admissible maps. Thus ${\Phi_X}(s)$ is a $Q$-topology on $X$.\\
\qquad \indent Let $\tau \in Q(X)$, ${\cal F}=\{f:(X, \tau)\rightarrow (S, u) | f$ is $Q$-continuous$\}$, and $s_\tau$ be the optimal lift of the family $\{f:X\rightarrow (S, u) | f\in{\cal F}\}$. Then clearly, $s_\tau \in C(X)$ and $f:(X, s_\tau) \rightarrow (S, u)$ is $C$-admissible for each $f\in{\cal F}$. This provides a function from $Q(X)$ to $C(X)$, which, it can be verified, is the inverse of $\Phi_X: C(X)\rightarrow Q(X)$. So, $C(X)$ and $Q(X)$ are in one-to-one correspondence. Also, condition (3) says that $f:(X, s)\rightarrow (Y, t)$ is $C$-admissible precisely when $q\circ f\in {\Phi_X}(s)$ for each $q\in {\Phi_Y}(t)$.

\end{document}